\documentclass[11pt,twoside,french]{article}

\usepackage{babel}
\usepackage{xypic}

\font\teneusm=eusm10 at 12pt
\font\seveneusm=eusm7 at 9pt
\font\fiveeusm=eusm5 at 7pt
\newfam\eusmfam
\textfont\eusmfam=\teneusm
\scriptfont\eusmfam=\seveneusm
\scriptscriptfont\eusmfam=\fiveeusm

\font\tenmsbm=msbm10 at 12pt
\font\sevenmsbm=msbm7 at 9pt
\font\fivemsbm=msbm5 at 7pt
\newfam\msbmfam
\textfont\msbmfam=\tenmsbm
\scriptfont\msbmfam=\sevenmsbm
\scriptscriptfont\msbmfam=\fivemsbm
\def\msbm{\tenmsbm\fam\msbmfam}

\newtheorem{proposition}[subsection]{Proposition}
\newtheorem{corollary}[subsection]{Corollaire}

\def\proofbox{\null\nobreak\hfill\rule{6pt}{6pt}}

\newenvironment{proof}[1]{\begin{trivlist} \item[] {\bf D\'emonstration
de~#1}}{\proofbox\end{trivlist}}

\newenvironment{numpar}{\begin{trivlist} \item[] {\addtocounter{subsection}{1}\bf \arabic{section}.\arabic{subsection}} }
{\end{trivlist}}

\newenvironment{encarttitre}[2]{%
  \begin{center}%
  {\bfseries #1\vspace{-.5em}}
  \end{center}%
  \small
  \begin{quotation} {\bf #2}}
 {\end{quotation}}

\newenvironment{encart}[1]{%
  \small
  \begin{quotation} {\bf #1}}
 {\end{quotation}}

\def\Pm{{\bf P}_m}
\def\P{{\bf P}}
\def\Pp{{{\bf P}'}}

\def\sF{{\cal F}}
\def\I{{\cal I}}

\def\piinv{\pi^{\ast}}
\def\pidir{\pi_{\ast}}
\def\pirec{\pi^{-1}}
\def\piadj{\pi^{!}}
\def\OP{{\cal O}_{\P}}
\def\OPm{{\cal O}_{\Pm}}
\def\OPp{{\cal O}_{\Pp}}
\def\OX{{\cal O}_{X}}
\def\OXp{{\cal O}_{X'}}

\def\OSp{{\cal O}_{S'}}

\def\canP{{\omega}_{\P}}

\def\canPp{{\omega}_{\Pp}}
\def\canX{{\omega}_{X}}
\def\canXp{{\omega}_{X'}}
\def\canS{{\omega}_{S}}

\def\canSp{{\omega}_{S'}}
\def\sHom{{\underline{\rm Hom}}}
\def\dim#1{{\rm dim}\left(#1\right)}
\def\codim#1{{\rm codim}\left(#1\right)}
\def\cH#1#2{{\rm H}^{#1}\!\left(#2\right)}
\def\ch#1#2{{\rm h}^{#1}\!\left(#2\right)}

\def\ms{\mu_{s}}
\def\ptt{\ {\rm pour\ tout}\ }

\def\map#1{\stackrel{#1}{\longrightarrow}}

\begin{document}
\pagestyle{myheadings}
\markboth{Sections hyperplanes et endomorphismes}{\rm G\'{e}om\'{e}trie alg\'{e}brique/\em Algebraic geometry}


\newpage
\null

\vskip 1em%
\begin{center}%
  \let \footnote \thanks
    {\LARGE Sections hyperplanes et endomorphismes de l'espace projectif \par}
    \vskip 1.5em%
    {\large
      \lineskip .5em%
      \begin{tabular}[t]{c}
        Guillaume Jamet\\
	{\small le 15 mars 1997}	
      \end{tabular}\par}%
\end{center}%
\par
\vskip 1.5em



\begin{encart}{R\'{e}sum\'{e}}
Nous montrons que toute section hyperplane d'une vari\'{e}t\'{e} image r\'{e}ciproque d'une
vari\'{e}t\'{e} lisse de dimension au moins \'{e}gale \`{a} 2 par un
endomorphisme (qui n'est pas un automorphisme) de l'espace projectif est
lin\'{e}airement compl\`{e}te. Ce r\'{e}sultat a un int\'{e}r\^{e}t
particulier dans le cadre des surfaces lisses de $\P_{4}$.
\end{encart}

\begin{encarttitre}{Hyperplane sections and endomorphisms of projective space}{Abstract}
We show that any hyperplane section of a variety which is the inverse
image of a smooth variety of dimension at least 2 by an
endomorphism (wich is not an automorphism) of the projective space, is
linearly complete. We stress the case of smooth surfaces in $\P_{4}$.
\end{encarttitre}


\section*{Introduction.}
G.~Ellingsrud et C.~Peskine ont d\'{e}montr\'{e} (\cite{GE-CP}) que,
except\'{e} pour un nombre fini de composantes du sch\'{e}ma de Hilbert, une
surface lisse de $\P_{4}$ est de type g\'{e}n\'{e}ral. D'autre part,
les seuls exemples de surfaces lisses irr\'{e}guli\`{e}res de $\P_{4}$
connus sont d'irr\'{e}gularit\'{e} \'{e}gale \`{a} 1 ou
2. Ces surfaces sont construites comme images r\'{e}ciproques par des
endomorphismes d'un nombre {\em fini} de surfaces irr\'{e}guli\`{e}res
connues. Et il est conjectur\'{e} de longue date que
l'irr\'{e}gularit\'{e} des surfaces de $\P_{4}$ est
born\'{e}e. Dans~\cite{CP}, Peskine conjecture que, except\'{e} pour un
nombre fini de composantes du sch\'{e}ma de Hilbert, la section
hyperplane g\'{e}n\'{e}rale d'une surface lisse de $\P_{4}$ est
lin\'{e}airement compl\`{e}te (cf. \S\ref{dec}). Il sugg\`{e}re aussi que la section
hyperplane  g\'{e}n\'{e}rale d'une surface obtenue comme image
r\'{e}ciproque par un endomorphisme de $\P_{4}$ d'une surface irr\'{e}guli\`{e}re,
est lin\'{e}airement compl\`{e}te, except\'{e} peut-\^{e}tre  pour un nombre fini de familles de surfaces
irr\'{e}guli\`{e}res. 
Nous montrons ici que toute section hyperplane d'une surface image
r\'{e}ciproque d'une surface par un
endomorphisme (qui n'est pas un automorphisme) est
lin\'{e}airement compl\`{e}te. Ce qui, vue la nature des surfaces
irr\'{e}guli\`{e}res connues, va dans le sens de la conjecture de
Peskine. La d\'{e}monstration repose sur une d\'{e}composition de
certaines images directes de fibr\'{e}s (\ref{decomp}) et conduit en
fait \`{a} un \'{e}nonc\'{e} g\'{e}n\'{e}ral sur les  sections
hyperplanes de vari\'{e}t\'{e}s images r\'{e}ciproques de
sous-vari\'{e}t\'{e}s lisses  par un endomorphisme de l'espace projectif (\ref{hyper2}). La
derni\`{e}re section regroupe quelques remarques sur l'adjonction de ces
vari\'{e}t\'{e}s.

\section{D\'{e}composition de $\pidir (\OXp (lH'))$.}
\label{dec}
Soit $n$ un entier $\geq 1$. Soient $\P$ et $\Pp$ des espaces projectifs
sur {\bf C} de dimension $n$ et $H$ (resp. $H'$) une section hyperplane
de $\P$ (resp. $\Pp$). Nous dirons qu'un sous-sch\'{e}ma $X$ d'un espace
projectif $\P$ est lin\'{e}airement complet si le morphisme de
restriction~: $\cH{0}{\P, \OP(H)}\map{}\cH{0}{X, \OX(H)}$ est surjectif.\\
Soient $\Pp \map{\pi} \P$ un morphisme d'image non r\'eduite \`a un point,
donc fini, surjectif et plat, $X \subset\P$ un sous-sch\'ema de
degr\'{e} $d$ et $X'=\pirec (X)\subset\Pp$ son image r\'eciproque. Si
$\piinv (\OP (H))=\OPp (k H')$, $X'$ est de degr\'e
$d.k^{\codim{X}}$. En particulier $\pi$ est de degr\'e $k^n$. 

Lorsque $X$ est localement  Cohen-Macaulay, on d\'{e}signera par
$\canX$ son faisceau dualisant, et, dans ce cas, $\pi$ \'{e}tant
fini et plat, $X'$ est aussi localement Cohen-Macaulay et son
faisceau dualisant est $\canXp=\piadj(\canX)$ (o\`{u} $\piadj(\sF)$ est
d\'{e}fini par $\pidir\piadj(\sF)=\sHom(\pidir(\OXp),\sF)$, cf.~\cite{RD}
ou~\cite[III.7]{RH}). On souhaite \'etudier les
faisceaux $\OXp (lH')$ et, lorsque $X$ est localement Cohen-Macaulay,
$\canXp (lH')$, et  en particulier leur cohomologie. La proposition suivante
et l'annulation des images directes sup\'{e}rieures de $\pi$ nous
permettent de <<travailler sur $X$>>.

\begin{proposition}\label{decomp}
Le fibr\'e $\pidir (\OXp (l H')),\ l\in {\msbm Z}$, se d\'ecompose
naturellement~:
$$
\pidir (\OXp (l H'))=\bigoplus_{d\in{\msbm Z}}E_{l,d}\otimes\OX (-d H)
$$ 
o\`u $E_{l,d}$ est le conoyau de la multiplication~:
$$
\cH{0}{\OP (H)}\otimes \cH{0}{\pidir (\OPp(l
H'))((d-1)H)}\map{}\cH{0}{\pidir (\OPp (l H'))(d H)},
$$
non nul exactement lorsque $-\left[ l/k\right]\leq d\leq\delta(n,k,l)$,
avec $\delta (n,k,l)=n+1+\left[ -\frac{n+1+l}{k}\right]$.
\end{proposition}
\begin{proof}{{\ref{decomp}}}
Comme $\pi$ est affine, $\pidir(\OXp (l H'))=\pidir(\OPp (l H'))_{\vert
X}$, on peut donc
supposer que $X=\P$. Or $\pidir(\OPp (l H'))$ est un fibr\'e sur $\P$ et
$$
\cH{i}{\pidir(\OPp (l H'))(d H)}=\cH{i}{\OPp ((l+k.d)H')}.
$$
Le terme de droite s'annule pour $0<i<n$ et tout entier $d$.
Donc $\pidir(\OPp (l H'))$ est une somme de fibr\'es en droites. Un tel
fibr\'{e} est engendr\'{e} en degr\'{e} $\leq\delta$ si et seulement si
$$
\cH{n}{\pidir (\OPp (l H'))((\delta-n) H)}=0,
$$
ce qui est \'{e}quivalent \`{a} $l+k(\delta-n)>-n-1$. Notons
$V_{l,d}=\cH{0}{\pidir(\OPp(lH'))(dH)}$ (la composante homog\`{e}ne de
degr\'{e} $l+kd$ de l'anneau gradu\'{e} de $\Pp$). Le diagramme commutatif suivant
(o\`{u} $l+kd\geq 0$ et les fl\`{e}ches horizontales
(resp. verticales) sont induites par la multiplication des sections sur
$\P$ (resp. $\Pp$))~:
$$
\xymatrix{
{\cH{0}{\OP(H)}\otimes V_{l,d-1}\otimes V_{0,1}}\ar[r]\ar@{->>}[d]
 & {V_{l,d}\otimes V_{0,1}}\ar[r]\ar@{->>}[d]
 & {E_{l,d}\otimes V_{0,1}}\ar[r]\ar@{->>}[d] & 0 \\
{\cH{0}{\OP(H)}\otimes V_{l,d}}\ar[r]
 & {V_{l,d+1}}\ar[r]
 & {E_{l,d+1}}\ar[r] & 0
}
$$ 
montre que $E_{l,d}=0$ entra\^{\i}ne $E_{l,d+1}=0$~. D'o\`{u} la d\'{e}composition annonc\'{e}e.
\end{proof}

\begin{numpar}\label{nondeg}
Lorsque $\pi$ n'est pas un automorphisme (i.e. $k>1$), on v\'{e}rifie
que $\cH{0}{\I_{X'}(H')}=0$, donc $X'$ est non d\'{e}g\'{e}n\'{e}r\'{e}.
\end{numpar}
\begin{numpar}\label{complet}
D'autre part si $\ch{1}{\I_{X}(d H)}=0\ptt d<s$ alors
$$
\ch{1}{\I_{X'}(l H')}=0\ptt l<sk.
$$
En effet, si $\sF$ est un faisceau coh\'{e}rent sur $X$, on a
$$\cH{i}{\piinv (\sF)(l H')}\simeq\cH{i}{\sF\otimes\pidir (\OXp (l H'))}.
$$
En particulier si $\ch{0}{\OX}=1$ et si $\pi$ n'est pas un automorphisme,
$X'$ est lin\'{e}airement complet. Remarquons cependant que $\ch{1}{\I_{X}(d H)}=0\ptt d$ si et seulement si $\ch{1}{\I_{X'}(l H)}=0\ptt l$.
\end{numpar}

\begin{numpar}
Si $X$ est lisse on sait d'apr\`{e}s le th\'{e}or\`{e}me de
transversalit\'{e} de Kleiman (voir~\cite{SK} ou~\cite[III.10.8, p. 273]{RH}) que, pour un
automorphisme lin\'{e}aire g\'{e}n\'{e}ral $\sigma$, $\pirec (\sigma X)$
est lisse.
Donc si $X$ est lisse et en position g\'{e}n\'{e}rale, $X'$ est lisse et
$\ch{i}{\OXp}=\ch{i}{\OX},\ i<\dim{X}$.
\end{numpar}
Dans toute la suite on supposera que $\pi$ n'est pas un
automorphisme (i.e. $k>1$).

\section{La section hyperplane.}

Soient $X\subset\P$ une vari\'{e}t\'{e} lisse connexe, $s$ une
section non-nulle de $\OPp(H')$ et $\ms$ la
multiplication par $s$~: $\OXp\map{\ms}\OXp(H')$.
\begin{proposition}\label{hyper1}
$\cH{i}{\ms}$ est injective pour $i<\dim{X}$.
\end{proposition}
Comme d'apr\`{e}s {1.2} et {1.3} l'image r\'{e}ciproque de $X$ est
lin\'{e}airement compl\`{e}te et non d\'{e}g\'{e}n\'{e}r\'{e}e, on en
d\'{e}duit~:
\begin{corollary}\label{hyper2}
Soient $X\subset\P$ une variet\'{e} lisse de dimension au moins \'{e}gale
\`{a} 2, $\pi :\Pp\map{}\P$ un endomorphisme qui n'est pas un
automorphisme  et $X'=\pirec(X)$. Alors pour tout hyperplan $H'$ coupant
$X'$ proprement, la section $X'\cap H'$ est lin\'{e}airement
compl\`{e}te. En particulier, si $X'$ est irr\'{e}ductible, toute section
hyperplane de $X$ est lin\'{e}airement compl\`{e}te.
\end{corollary}
\begin{proof}{\ref{hyper1}}
D'apr\`{e}s~\ref{decomp}, 
$\pidir(\ms)$ se d\'{e}compose sous la forme: $\pidir(\ms)=\left(\ms
^{p,q}\right)_{p,q\geq 0}$ avec~:
$$
\ms ^{p,q} : E_{0,q}\otimes\OX (-q H)\map{}E_{1,p}\otimes\OX (-p H)
$$
Remarquons que $\ms ^{0,0}$ s'identifie \`{a}~:
$$
\OX\map{s\otimes 1}\cH{0}{\OPp(H')}\otimes\OX
$$
En appliquant le th\'{e}or\`{e}me d'annulation de Kodaira on
voit que, pour $i<\dim{X}$, $\cH{i}{\pidir(\ms)}$ se r\'{e}duit \`{a}~:
$$
\cH{i}{\OX}\map{s\otimes 1}\cH{0}{\OPp(H')}\otimes\cH{i}{\OX}
$$
qui est \'{e}videmment injective.
\end{proof}
Sans supposer $X$ lisse, on d\'{e}montre de
la m\^{e}me mani\`{e}re~:
\begin{proposition}
Soit $j$ un entier. Si $\cH{i}{\OX(-dH)}=0$ pour tout entier $i<j$ et tout
entier $d$ tel que $0<d\leq\delta(n,k,0)$ alors
$\cH{i}{\ms}$ est injective pour tout entier $i<j$.
\end{proposition}

\section{L'application canonique.}

On suppose toujours que $X$ est lisse, connexe et en position
g\'{e}n\'{e}rale. Si $0\leq l<k$ et $\delta_{l}=\delta(n,k,l)$, on a $0<\delta_{l} < n+1$.
D'apr\`{e}s~\ref{decomp}~:
$$
\pidir(\canXp(-lH'))=\sHom(\pidir(\OXp(lH')),\canX)=(E^{\vee}_{l,0}\otimes\canX)\oplus\ldots\oplus(E^{\vee}_{l,\delta_{l}}\otimes\canX(\delta_{l}H))
$$
Un \'{e}l\'{e}ment non nul de $E^{\vee}_{l,\delta_{l}}$ d\'{e}finit
naturellement un morphisme non nul~: $$\piinv(\canX(\delta_{l}
H))\map{}\canXp(-lH')$$
\begin{numpar}
Si $l=0$, on obtient une factorisation~:
$$
\xymatrix{
{\P(\cH{0}{\canXp})}\ar@{-->}[d]
 & X'\ar@{-->}[l]\ar[d] \\
{\P(\cH{0}{\canX(\delta_{0} H)})}
 & X\ar@{-->}[l]
}
$$ 
La fl\`{e}che verticale de gauche est une projection.
\end{numpar}
\begin{numpar}
Si $l=1$ et $\delta_{1}\geq\dim{X}$, on sait (cf.~\cite[8.8.5, p. 239]{BS}) que
$\canX(\dim{X}H)$ est  engendr\'{e} par ses sections sauf si
$(X,\OX(H))=(\Pm,\OPm(1))$. Donc, hormis ce cas, $\cH{0}{\canXp(-H')}$
est $\neq 0$. Si $s$ est une section non nulle de $\canXp(-H')$, le
syst\`{e}me lin\'{e}aire $s\otimes\cH{0}{\OXp(H')}\subset\cH{0}{\canXp}$,
d\'{e}fini une application rationnelle qui est un plongement hors du
lieu des z\'{e}ros de $s$ et donc l'application canonique de $X'$ est birationnelle.
\end{numpar}
En particulier~:
\begin{proposition}\label{typegen}
Soient $S\subset\P_{4}$ une surface lisse et $\pi :\P_{4}\map{}\P_{4}$ un
endomorphisme qui n'est pas un automorphisme. Alors, si $S$ est en
position g\'{e}n\'{e}rale, $S'=\pirec(S)$ est lisse et son fibr\'{e}
canonique est tr\`{e}s ample, sauf si $S$ est un plan et $k=2$.
\end{proposition}
\begin{numpar}
La surface  $S'$ est donc de type g\'{e}n\'{e}ral,
sauf dans le cas d'exception de la proposition, $S'$ \'{e}tant alors une
surface de Del Pezzo.
\end{numpar}
\begin{proof}{\ref{typegen}}
Comme $N_{S'}\Pp=\piinv(N_{S}\P)$, alors 
$\canSp\otimes\piinv(\canSp)^{\vee}=(\canPp\otimes\piinv(\canP)^{\vee})_{\vert
S'}=\OSp(5(k-1)H')$. Donc
$\canSp=\piinv(\canS(3H))\otimes\OSp((2k-5)H')$. Cela d\'{e}montre le
r\'{e}sultat puisque, sauf si $S$ est un plan, $\canS(2H)$ est
engendr\'{e} par ses sections.
\end{proof}

\vfil

\hfill\vbox{
\hbox{Guillaume Jamet}\smallskip
\hbox{Universit\'{e} Pierre et Marie Curie}
\hbox{Institut de Math\'{e}matiques}
\hbox{Analyse Alg\'{e}brique, Case 82}
\hbox{Tour 46--00, $3^{e}$ \'{e}tage}
\hbox{4, place Jussieu}
\hbox{F-75252 PARIS Cedex 05}\smallskip
\hbox{T\'{e}l~: 01 44 27 26 90}
\hbox{Fax~: 01 44 27 26 87}
\hbox{e-mail~: {\tt jamet@math.jussieu.fr}}
}
  
\end{document}